\documentclass[12pt,singlespace,oneside]{article}
\usepackage{amsmath}
\usepackage{latexsym}

\input{tcilatex}
\begin{document}

\title{Helmholtz Theorem for Differential Forms in 3-D Euclidean Space%
\thanks{%
To Dr. Howard Brandt, for his contribution to the improvement of my papers
and books.}}
\author{Jose G. Vargas\thanks{%
PST Associates, 138 Promontory Rd, Columbia, SC 29209-1244. USA,\newline
\hspace*{0.5cm} josegvargas@earthlink.net} \ }
\date{April 10, 2014.}
\maketitle

\begin{abstract}
There are significant differences between Helmholtz and Hodge's
decomposition theorems, but both share a common flavor. This paper is a
first step to bring them together.

We here produce Helmholtz theorems for differential $1-$forms and $2-$forms
in 3-D Euclidean space, $E_{3}.$ We emphasize their common structure in
order to facilitate the understanding of another paper, soon to be made
public, where a Helmholtz theorem for arbitrary differential forms in
arbitrary Euclidean space is presented and which allows one to connect
(actually to derive from it) an improvement of Hodge's decomposition theorem.
\end{abstract}

\section{Context of this Paper}

This paper is a first step in connecting the Helmholtz and Hodge
decomposition theorems. The first of these pertains to vector fields
(objects of grade one) in 3-D Euclidean space, $E_{3}.$ It involves the
action of the del operator through vector product. This product is specific
to three dimensions. It does not exist in arbitrary dimension. Hence the
present Helmholtz theorem is a peculiarity of 3-D.

The second theorem pertains to the vast subject of arbitrary differential
forms (in compact oriented Riemannian manifolds) of arbitrary dimension. The
first one is of interest mainly for physicists in general, and can be found
in the first chapter of a book for them on mathematical methods \cite%
{ArfkenW}. The second one is of interest primordially for mathematicians and
mathematical physicists, specially those concerned with any of the
overlapping subjects of de Rham\cite{de Rham}-Hodge\cite{Hodge} theory,
elliptic operators, cohomology, etc. They appear in the thick of books where
at last one of those subjects is considered.

This author has made progress in relating those two theorems. The cuspid of
his results is constituted by the following:

(a) A decomposition theorem for differential forms in Euclidean spaces, $%
E_{n}$, of arbitrary dimension in \textit{Helmholtz format}. This
terminology is justified for making a difference with the \textit{Hodge
format}, where there are only two terms in the decomposition (a
would-be-third term goes to zero if and when relevant quantities decay
sufficiently fast at infinity). That terminology is also justified because,
as in the original theorem and unlike in Hodge's case, the terms of the
decomposition are specified.

(b) An application of (a) through the use of the embedding theorem of Schl%
\"{a}fli\cite{Schlafli} -Janet\cite{Janet} -Cartan\cite{Cartan} to improve
on the original Hodge theorem \cite{Hodge}. \ For modern treatments, see %
\cite{Choquet}, \cite{Morita}, \cite{Frankel} and \cite{Warner}. Let us
recall that the original theorem states that a differential form in an
oriented compact Riemannian manifolds can be decomposed into three
components that are respectively closed, co-closed and harmonic, but it does
not state what those components are for given differential form. Our version
of the theorem does. For the proof, we resort to K\"{a}hler's calculus (KC) %
\cite{K60}, \cite{K62}, which generalizes Cartan's.

None of these results is presented in this paper (I), but in a paper
intended to follow (II), which requires more mathematical sophistication.
But I may have enough contents to satisfy the interest of most physicists. A
second reason is that this paper seeks to pique the interest of a potential
sponsor of II for the mathematics-analysis section of the arXiv, where such
paper would presumably be assigned if sponsored. It is important to realize
that an important indirect goal of both papers is to call attention to the
fact that KC\ is a formidable tool that is being overlooked. It provides
opportunities to obtain substantial results in dealing with issues at which
one need not be an expert. That is precisely what makes it so important.

\section{Introduction}

Cartan's calculus enriched by the Hodge duality operation may suffice to
follow this paper, but KC provides a much richer context for the arguments
and computations. Consider, for instance, the co-differential. In the so
enriched Cartan calculus, this concept involves the metric through Hodge
duality. In KC, it is a more general concept introduced through the
connection. Both definitions coincide for the Levi-Civita connection, but
not in general.

The present author is not aware of whether every particular result he uses
is known by practitioners of Cartan's calculus enriched by the Hodge duality
operation. Hence, in section 3, we have introduced well known results (and
some which may not be as well known) from a KC perspective. Our source is
any of the papers \cite{K60} and \cite{K62}. Their author owes much to De
Rham and Hodge, as he recognized in the first paragraph of \cite{K60}.

In section 4, we speak of the Helmholtz format. These considerations, later
enriched as we make progress in the paper, will allow us to see very early
what the Helmholtz format is for general differential forms in euclidean
spaces.

In section 5, we derive Helmholtz' theorem for differential $2-$forms. The
proof mimics the one for vector fields. In section 6, we derive Helmholtz'
theorem for differential $2-$forms, using the fact that any such
differential form has a differential $1-$form as dual. Hence, from the
theorem for the last ones, we obtain the theorem for the first ones. This is
instructive in the following specific regard.

Recall that, in the derivation of Helmholtz theorem for vector fields, the
treatment of one of the two terms is more complicated than the derivation of
the other term. This asymmetry carries to the corresponding theorem for
differential $1-$forms. The complication is reversed for differential $2-$%
forms. The reason behind it is that the co-differential of a $1-$form is a $%
0-$form, and the exterior differential of a $2-$form is a $3-$form, which,
in dimension 3, is almost like a $0-$form for certain purposes. In higher
dimension both terms of the decomposition will present similar difficulty.
Which brings us already to state (without proof) the following qualitative
feature of the improvement of our Hodge theorem. The harmonic term is in
turn the sum of two terms, each of them harmonic and respectively related to
the closed and co-closed terms.

Finally, in section 7, we anticipate how the generalization of Helmholtz
theorem of which we spoke above will look like.

How to read this paper? It will be a function of a reader's knowledge. If
familiar with KC, one can comfortably jump to section 5. If not, it may help
to read every section.

\section{Calculus of Differential Forms}

Our differential forms are integrands (known as currents in some of the
literature), not skew-symmetric multilinear functions of vectors. Stokes
generalized theorem is directly about integrals, and only indirectly about
skew-symmetric multilinear functions of vectors. Its underlying algebra,
when dealing only scalar-valued differential forms, is Clifford algebra, but
is known as K\"{a}hler's algebra when it refers to differential forms. Thus,
not all calculi based on Clifford algebra are equivalent. Alternatives are
the calculi by Dirac \cite{Dirac} and by Hestenes \cite{Hestenes} , which
are both based on tangent algebra. Worth mentioning is that the problem with
negative energies does not arise in KC \cite{K62}.

From what has been said, K\"{a}hler's algebra is built upon the module of
differential forms spanned by $(dx^{i})$ and defined by 
\begin{equation}
dx^{i}dx^{j}+dx^{j}dx^{i}=2g^{ij}.  \label{1}
\end{equation}%
When at least one of two factors in a Clifford product is a differential
1-form, $\alpha $, the identity 
\begin{equation}
\alpha \Gamma \equiv \frac{1}{2}(\alpha \Gamma +\Gamma \alpha )+\frac{1}{2}%
(\alpha \Gamma -\Gamma \alpha )  \label{2}
\end{equation}%
allows one to define the two terms on the right of (2) as the interior and
exterior product. Which one is which depends on the grade of $\Gamma $. If $%
\Gamma $ is another 1-form, $\beta $, we have 
\begin{equation}
\alpha \beta =\alpha \wedge \beta +\alpha \wedge \beta   \label{3}
\end{equation}%
where 
\begin{subequations}
\begin{equation}
\alpha \wedge \beta \equiv \frac{1}{2}(\alpha \beta -\beta \alpha ),
\label{4a}
\end{equation}%
\begin{equation}
\alpha \cdot \beta \equiv \frac{1}{2}(\alpha \beta +\beta \alpha ).
\label{4b}
\end{equation}%
Those specific ``$\wedge $'' and ``$\cdot $'' products are of respective
grades two and zero. These equations apply in particular to when $\alpha $
and $\beta $ are the differentials of the coordinate functions.

Let $w$ denote the unit differential of highest possible grade. In Cartesian
coordinates in 3-D Euclidean space, $E_{3}$, it can be given as 
\end{subequations}
\begin{equation}
w=dxdydz.  \label{5}
\end{equation}%
Since the exterior product raises the grade, the product ``$w\wedge ...$''
by the algebra (except scalars) yields zero. Thus, we have%
\begin{equation}
\Gamma w=\Gamma \wedge w+\Gamma \cdot w=\Gamma \cdot w,  \label{6}
\end{equation}%
for non-scalars. Multiplication by $w$ is called the obtaining of the Hodge
dual, which corresponds in the tensor calculus to contracting with the
Levi-Civita tensor. In $E_{3}$, $w$ commutes with the whole algebra for that
space, and $w^{2}=-1.$

The K\"{a}hler operator, which we shall represent with the symbol $\partial $%
, is the sum of a part, $d,$ which raises the grade by one, and a part, $%
\delta $, which lowers the grade by one:%
\begin{equation}
\partial =d+\delta .  \label{7}
\end{equation}%
Readers may refer to it as Dirac's operator, but we prefer to keep our
distance from anything Dirac's since it may occasionally induce one into
error due to different contexts.

When the connection of the manifold is the Levi-Civita connection (LCC), $%
\delta $ of a scalar-valued differential form is called co-differential. The
operation $d$ is exterior differentiation ---exterior covariant
differentiation if applied to tensor-valued or Clifford-valued differential
forms. Neither $\partial $ nor $\delta $ satisfy the standard Leibniz rule.
And yet, K\"{a}hler refers to $\partial $ (for which he uses the symbol $%
\delta $) as interior differentiation. There is strong reason for using the
term differentiation in spite of the issue about the Leibniz rule. But K\"{a}%
hler's better theory should not always be constrained by how terms are used
in other theories.

Under the LCC, differential forms $\Gamma $ of arbitrary grade satisfy%
\begin{equation}
\delta \Gamma =\left( -1\right) ^{\tbinom{n}{2}}d(\Gamma w)w.  \label{8}
\end{equation}%
In dimension three, we readily get%
\begin{equation}
\delta \Gamma =-wd(\Gamma w),  \label{9}
\end{equation}%
which is not a definition but a theorem in KC. In arbitrary dimension,%
\begin{equation}
(\delta \Gamma )w=d(\Gamma w)  \label{10}
\end{equation}%
since $w^{2}=\left( -1\right) ^{\tbinom{n}{2}}.$

The Laplacian of differential forms is defined as 
\begin{equation}
\partial \partial \Gamma =(d+\delta )(d+\delta )\Gamma =(dd+d\delta +\delta
d+\delta \delta )\Gamma .  \label{11}
\end{equation}%
For scalar-valued differential forms, $dd\Gamma $ is zero. If, in addition,
the connection is Levi-Civita's, $\delta \delta \Gamma =0.$ Thus, 
\begin{equation}
\partial \partial \Gamma =d\delta \Gamma +\delta d\Gamma ,
\end{equation}%
and further, 
\begin{equation}
\partial \partial f=\delta df
\end{equation}%
for 0-forms, $f$.

All the above is well known. What follows may rarely be known. K\"{a}hler
defines covariant derivatives of ($p,q$)-valued differential $r-$forms,
inhomogeneous in general. Their components have three series of indices, $p$%
, $q$ and $r$, the $q$ and $r$ series being of subscripts; $q$ is for
multilinear functions of vectors and $r$ for functions of hypersurfaces,
i.e. integrands. Thus curvature is (1,1)-tensor valued differential 2-form.
As a differential 2-form, it is a function of surfaces, surfaces defined by
pairs of curves with the same origin and end along which we transport a
vector. The second number one in (1,1) refers to the vector being
transported. Curvature is then a function of a pair of vector field and
surface. This function is vector-valued, meaning the following. It is
evaluated on the vector field and then evaluated (read integrated) on the
surface. As already mentioned by Cartan, integration surfaces must be
infinitesimal, unless the affine curvature is zero.

We only need here scalar-valuedness ($p=$ $q=0$). The definition of
covariant derivative then simply is 
\begin{equation}
d_{h}u\equiv \frac{\partial u}{\partial x^{h}}-\omega _{h}^{k}\wedge e_{k}u,
\label{14}
\end{equation}%
where $e_{k}u$ is $\omega _{k}\cdot u$. K\"{a}hler then defines $\partial u$
as%
\begin{equation}
\partial u\equiv dx^{h}\vee d_{h}u=du+\delta u  \label{15}
\end{equation}%
where%
\begin{equation}
du\equiv dx^{h}\wedge d_{h}u\text{, \ \ \ \ \ \ \ }\delta u\equiv
dx^{h}\cdot d_{h}u.  \label{16}
\end{equation}%
Consider Euclidean space. It is clear that $du=dx^{h}(\partial u/\partial
x^{h})$ since $dx^{h}\wedge \omega _{h}^{k}=ddx^{k}=0$. Notice that, in
Cartesian coordinates (which presupposes Euclidean space), $\delta
u=dx^{h}\cdot (\partial u/\partial x^{h})$, since $\omega _{h}^{k}$ is then $%
0.$ If $u$ were a differential $1-$form, the resulting expression would be
the same as the divergence of a vector field that, in terms of a constant
orthonormal frame field, has the same components as the differential form in
Cartesian coordinates.

K\"{a}hler proves that $\delta u$ so defined coincides with the
co-differential. $d_{h}u$ satisfies the Leibniz rule, but $\partial $ and $%
\delta $ do not. They rather satisfy%
\begin{equation}
\partial (u\wedge v)=\partial u\wedge v+\eta u\wedge \partial v+e^{h}u\wedge
d_{h}v+\eta d_{h}u\wedge e^{h}v,  \label{17}
\end{equation}%
and%
\begin{equation}
\delta (u\wedge v)=\delta u\wedge v+\eta u\wedge \delta v+e^{h}u\wedge
d_{h}v+\eta d_{h}u\wedge e^{h}v.  \label{18}
\end{equation}%
In view, however, of how $\partial $, $d$ and $\delta $ emerge from the
covariant derivative, it seems natural to refer to all of them as
differentiations, specially since we also need a name for $\delta u$
whenever the connection is not the LCC. Hence, we shall refer to $\partial $%
, $d$ and $\delta $ as, respectively, the K\"{a}hler, exterior and interior
differentials.

A differential form such that $d_{h}u$ is zero is called constant
differential. It will be denoted as $c$. They have the property that 
\begin{equation}
\partial (\Gamma c)=(\partial \Gamma )c,  \label{19}
\end{equation}%
and, in particular 
\begin{equation}
\partial (\Gamma w)=(\partial \Gamma )w.  \label{20}
\end{equation}%
Of great importance is that all polynomials in $(dx,dy,dz)$ with constant
coefficients are constant differentials, examples being $w$ and the $dx,dy$
and $dz$ themselves.

\section{Perspective on Helmholtz Theorem}

Helmholtz theorem states that smooth vector fields decaying sufficiently
fast at infinity can be written as a sum%
\begin{equation}
\boldsymbol{v}=-\boldsymbol{\nabla}\phi +\boldsymbol{\nabla}\times %
\boldsymbol{A}  \label{21}
\end{equation}%
with 
\begin{subequations}
\begin{equation}
\phi (\boldsymbol{r})=\frac{1}{4\pi }\int \frac{\boldsymbol{\nabla}^{\prime
}\cdot \boldsymbol{v}(\boldsymbol{r}^{\prime })}{r_{12}}dV^{\prime }
\label{22a}
\end{equation}%
\begin{equation}
\boldsymbol{A}(\boldsymbol{r})=\frac{1}{4\pi }\int \frac{\mathbf{\nabla }%
^{\prime }\times \boldsymbol{v}(\boldsymbol{r}^{\prime })}{r_{12}}dV^{\prime
}  \label{22b}
\end{equation}%
where $\nabla ^{\prime }$ refers to differentiation with respect to primed
coordinates, where $V^{\prime }$ is the volume element in primed coordinates
and where 
\end{subequations}
\begin{equation}
r_{12}\equiv \lbrack (x-x^{\prime })^{2}+(y-y^{\prime })^{2}+(z-z^{\prime
})^{2}]^{1/2}.  \label{23}
\end{equation}

We shall denote $\boldsymbol{v}(\boldsymbol{r}^{\prime })$ as $\boldsymbol{v}%
^{\prime }$. If we substitute (22) in (21), we get 
\begin{equation}
\boldsymbol{v}=\frac{1}{4\pi }\left[ -\boldsymbol{\nabla}\int \frac{1}{r_{12}%
}(\mathbf{\nabla }^{\prime }\cdot \boldsymbol{v}^{\prime })dV^{\prime }+%
\boldsymbol{\nabla}\times \int \frac{1}{r_{12}}\boldsymbol{\nabla}^{\prime
}\times \boldsymbol{v}^{\prime }dV^{\prime }\right] ,
\end{equation}%
whose format vis-a-vis the del is 
\begin{equation}
-\mbox{grad}\ldots \mbox{div}^{\prime }\boldsymbol{v}^{\prime }+\mbox{curl}%
\ldots \mbox{curl}^{\prime }\boldsymbol{v}^{\prime }.  \label{25}
\end{equation}%
Recalling that the Laplacian, $\Delta $, satisfies 
\begin{equation}
-\Delta =-\mbox{grad}\text{ }\mbox{div}+\mbox{curl}\text{ }\mbox{curl},
\label{26}
\end{equation}%
we are led to consider 
\begin{equation}
\partial \partial =-d\delta \alpha -\delta d\alpha   \label{27}
\end{equation}%
as replacement for (26) [strict parallelism between differentiations of
vector fields and differential $1-$forms would lead us to blindly write $%
-\Delta \boldsymbol{v}$ equal to $-\mbox{grad}$ $\mbox{div}\boldsymbol{v}-%
\mbox{div}$ $\mbox{grad}\boldsymbol{v},$ which is nonsense since vector
fields do not have gradients]. We thus consider the format 
\begin{equation}
-d\ldots \delta ^{\prime }\alpha ^{\prime }-\delta \ldots d^{\prime }\alpha
^{\prime }  \label{28}
\end{equation}%
for Helmholtz theorem for differential 1-forms, $\alpha $' being to $\alpha $
what $\boldsymbol{v}(\boldsymbol{r}^{\prime })$ is to $\boldsymbol{v}(%
\boldsymbol{r})$.

Those observations could be used to make the corresponding changes in (24)
to obtain Helmholtz theorem for differential 1-forms. But, in order to get
confidence with computations with them beyond the most trivial, we shall
formulate the theorem and proceed to prove it as one does in the vector
calculus. We shall thus need a uniqueness theorem like the one according to
which a vector field (respectively a differential 1-form, $\alpha $) is
uniquely defined in a region by its divergence and curl and its normal
component over the boundary (respectively its $\delta $ and $d$
''derivatives'' and its components at the boundary). The proof resorts to
showing that the difference between two hypothetical solutions, $\alpha _{1}$
and $\alpha _{2}$, have zero $d$ and $\delta $ derivatives. By the first of
these annulments $\alpha _{1}-\alpha _{2}$ is closed and, therefore, locally
(meaning ``not necessarily globally'') exact, i.e. $\alpha _{1}-\alpha
_{2}=df.$ We then have 
\begin{equation}
0=\delta (\alpha _{1}-\alpha _{2})=\delta df=\partial ^{2}f.  \label{29}
\end{equation}

We are now at a point similar to when, in the proof of the uniqueness
theorem in the vector calculus, one resorts to Green's theorem. The namesake
theorem in KC's is far more comprehensive, but it implies in particular that 
\begin{equation}
\int_{R}w(f\partial ^{2}g+\partial f\cdot \partial g)=\int_{\partial
R}(fdg)w.  \label{30}
\end{equation}%
We specialize this equation to $f=g$. We assume vanishing $f$ (i.e. $\alpha
_{1}-\alpha _{2}=0$) at the boundary. Equation (30) then yields, using (29),%
\begin{equation}
0=(\partial f)^{2}=\partial f=df,  \label{31}
\end{equation}%
all over the region considered. Hence $\alpha _{1}=\alpha _{2}$. The
uniqueness theorem has been proved.

\section{Helmholtz Theorem for Differential 1-Forms in 3-D Euclidean Space}

Helmholtz theorem: In Cartesian coordinates in $E_{3}$, differential $1-$%
forms that are smooth and vanish sufficiently fast at infinity can be
written as 
\begin{equation}
\alpha (\boldsymbol{r})=-\frac{1}{4\pi }dI^{0}-\frac{1}{4\pi }\delta
(dx^{j}dx^{k}I^{i}),  \label{32}
\end{equation}%
\begin{equation}
I^{0}\equiv \int \frac{1}{r_{12}}(\delta ^{\prime }\alpha ^{\prime
})w^{\prime },\ \ \ \ I^{i}\equiv \int \frac{1}{r_{12}}d^{\prime }\alpha
^{\prime }\wedge dx^{\prime i},  \label{33}
\end{equation}%
with summation over the three cyclic permutations of 1,2,3.

\textbf{Proof}: By the uniqueness theorem and the annulment of $dd$ and $%
\delta \delta $, the proof reduces to showing that $\delta $ and $d$ of
respective first and second terms on the right hand side of (32) yield $%
d\alpha $ and $\delta \alpha .$

Since $\delta dI^{0}=\partial \partial I^{0}$, we write $-(1/4\pi )\delta
dI^{0}$ as 
\begin{equation}
\frac{-1}{4\pi }\partial \partial I^{0}=\frac{-1}{4\pi }\int_{E_{3}^{\prime
}}(\partial \partial \frac{1}{r_{12}})(\delta ^{\prime }\alpha ^{\prime
})w^{\prime }=\frac{-1}{4\pi }\int_{E_{3}^{\prime }}(\partial ^{\prime
}\partial ^{\prime }\frac{1}{r_{12}})(\delta ^{\prime }\alpha ^{\prime
})w^{\prime }=\delta \alpha ,  \label{34}
\end{equation}%
after using the relation of $\partial \partial $ to the Dirac distribution.

For the second term, we use that $d\delta =\partial \partial -\delta d$ when
acting on $dx^{j}dx^{k}I^{i}.$ We move $\partial \partial $ past $%
dx^{j}dx^{k}$. Let $\alpha $ be given as $a_{l}(x)dx^{l}$ in terms of the
same coordinate system. We get $d^{\prime }\alpha ^{\prime }\wedge
dx^{\prime i}=(a_{k}^{\prime },_{j}-a_{j}^{\prime },_{k})w^{\prime }$. The
same property of $\partial \partial $ now allows us to obtain $d\alpha .$

For the second part of the second term, we apply $\delta d$ to $%
dx^{j}dx^{k}I^{i}:$%
\begin{equation*}
\delta d(dx^{j}dx^{k}I^{i})=\delta \left( w\frac{\partial I^{i}}{\partial
x^{i}}\right) =wd\left( \frac{\partial I^{i}}{\partial x^{i}}\right) =wdx^{l}%
\frac{\partial ^{2}I^{i}}{\partial x^{i}\partial x^{l}}=
\end{equation*}%
\begin{equation}
=wdx^{l}\int_{E_{3}^{\prime }}\left[ \frac{\partial ^{2}}{\partial x^{\prime
i}\partial x^{\prime l}}\left( \frac{1}{r_{12}}\right) \right]
(a_{k}^{\prime },_{j}-a_{j}^{\prime },_{k})w^{\prime }.  \label{35}
\end{equation}%
We integrate by parts with respect to $x^{\prime i}.$ One of the two
resulting terms is:%
\begin{equation}
wdx^{l}\int_{E_{3}^{\prime }}\frac{\partial }{\partial x^{\prime i}}\left[ 
\frac{\partial \left( \frac{1}{r_{12}}\right) }{\partial x^{\prime l}}%
(a_{k}^{\prime },_{j}-a_{j}^{\prime },_{k})\right] w^{\prime }.  \label{36}
\end{equation}%
Application to this of Stokes theorem yields%
\begin{equation}
wdx^{l}\int_{\partial E_{3}^{\prime }}\frac{\partial \left( \frac{1}{r_{12}}%
\right) }{\partial x^{\prime l}}(a_{k}^{\prime },_{j}-a_{j}^{\prime
},_{k})dx^{\prime j}dx^{\prime k}.  \label{37}
\end{equation}%
It vanishes for sufficiently fast decay at infinity.

The other term resulting from the integration by parts is%
\begin{equation}
-wdx^{l}\int_{E_{3}^{\prime }}\frac{\partial \left( \frac{1}{r_{12}}\right) 
}{\partial x^{\prime l}}\frac{\partial }{\partial x^{\prime i}}%
(a_{k}^{\prime },_{j}-a_{j}^{\prime },_{k})w^{\prime },  \label{38}
\end{equation}%
which vanishes identically (perform the $\frac{\partial }{\partial x^{\prime
i}}$ differentiation and sum over cyclic permutations). The theorem has been
proved.

\section{Helmholtz theorem for differential $2-$forms in $E_{3}$}

The theorem obtained for differential $1-$forms, here denoted as $\alpha $,
can be adapted to differential $2-$forms, $\beta $, by defining $\alpha $
for given $\beta $ as%
\begin{equation}
\alpha \equiv w\beta \text{, \ \ \ \ \ }\beta =-w\alpha .  \label{39}
\end{equation}%
Then, clearly,%
\begin{equation}
w\delta (w\beta )=-d\beta \text{, \ \ \ \ \ \ }wd\beta =\delta (w\beta ).
\label{40}
\end{equation}%
Helmholtz theorem for differential $1-$forms can then be written as%
\begin{equation}
w\beta =-\frac{1}{4\pi }d\left( \int_{E_{3}}\frac{\delta ^{\prime
}(w^{\prime }\beta ^{\prime })}{r_{12}}w^{\prime }\right) -\frac{1}{4\pi }%
\delta \left( dx^{jk}\int_{E_{3}}\frac{d^{\prime }(w^{\prime }\beta ^{\prime
})\wedge dx^{\prime i}}{r_{12}}\right) ,  \label{41}
\end{equation}%
and, therefore,%
\begin{equation}
\beta =\frac{1}{4\pi }wd\left( \int_{E_{3}}\frac{\delta ^{\prime }(w^{\prime
}\beta ^{\prime })}{r_{12}}w^{\prime }\right) +\frac{1}{4\pi }w\delta \left(
dx^{jk}\int_{E_{3}}\frac{d^{\prime }(w^{\prime }\beta ^{\prime })\wedge
dx^{\prime i}}{r_{12}}\right) .  \label{42}
\end{equation}%
The integrals are scalar functions of coordinates $x$. We shall use the
symbol $f$ to refer to them in any specific calculation. In this way, steps
taken are more easily identified.

The first term in the decomposition of $\beta $, we transform as follows:%
\begin{equation}
wdf=(\partial f)w=\partial (fw)=\delta (fw),  \label{43}
\end{equation}%
where we have used that $w$ is a constant differential.

For the second term, we have:%
\begin{equation}
w\delta (dx^{jk}f)=w\partial \lbrack wdx^{i}f)]-wd[fdx^{jk})].  \label{44}
\end{equation}%
The first term on the right is further transformed as%
\begin{equation}
w\partial (wdx^{i}f)=wwdx^{i}\partial f=-dx^{i}df,  \label{45}
\end{equation}%
where we have used that $wdx^{i}$ is a constant differential, which can be
taken out of the $\partial $ differentiation. For the other term, we have%
\begin{equation}
-wd(fdx^{jk})]=-wdf\wedge dx^{jk}=-wf,_{i}w=f,_{i}=dx^{i}\cdot df.
\label{46}
\end{equation}%
From the last three equations, we get%
\begin{equation}
w\delta (dx^{jk}f)=-dx^{i}df+dx^{i}\cdot df=-dx^{i}\wedge df=d(dx^{i}f).
\label{47}
\end{equation}

In order to complete the computation, we have to show that $d(w\beta )\wedge
dx^{i}$ can be written as $\delta \beta \wedge dx^{jk}.$ This can be shown
easily by direct calculation. Let $\alpha $ be given as $a_{i}dx^{i}.$ Then $%
d(w\beta )\wedge dx^{1}=d\alpha \wedge dx^{1}=(a_{3,2}-a_{2,3})w.$ On the
other hand, $\beta =-a_{i}dx^{jk}$ and%
\begin{equation}
\delta \beta =(a_{3,2}-a_{2,3})dx^{1}+cyclic\text{ }permutations.
\end{equation}%
Hence $\delta \beta \wedge dx^{23}=(a_{3,2}-a_{2,3})w$ and, therefore,%
\begin{equation}
d(w\beta )\wedge dx^{1}=d\alpha \wedge dx^{1}=\delta \beta \wedge dx^{23},
\end{equation}%
and similarly for the cyclic permutations of the indices.

\section{Anticipation of Helmholtz theorem for differential $r-$forms in $%
E_{n}$}

The discussion of section 4 together with the development of the proofs of
sections 5 and 6 makes it obvious what Helmholtz theorem for differential $r-
$forms will look like, namely

\begin{equation}
\alpha _{r}(x)=\mu d\left[ dx^{i_{1}...i_{r-1}}I^{\delta }\right] +\mu
\delta \left[ dx^{k_{1}...k_{r+1}}I^{d}\right] ,
\end{equation}%
\begin{equation}
I^{\delta }\equiv \int_{E_{n}^{\prime }}\frac{1}{r_{12}^{\lambda }}(\delta
^{\prime }\alpha _{r}^{\prime })\wedge dx^{\prime j_{1}...j_{n-r+1}},
\end{equation}%
\begin{equation}
I^{d}\equiv \int_{E_{n}^{\prime }}\frac{1}{r_{12}^{\lambda }}(d^{\prime
}\alpha _{r}^{\prime })\wedge dx^{\prime l_{1}...l_{n-r+1}},
\end{equation}%
$n$ being the dimension of the Euclidean space where $\alpha _{r}(x)$ would
have been defined, or the dimension of a still larger Euclidean space. The
integrations are performed over the chosen $E_{n}$ space.  Summation over
the basis made by $\frac{n!}{(r-1)!(n-r+1)!}$ independent basis elements $%
dx^{i_{1}\text{...}i_{r-1}}$, and the basis made by $\frac{n!}{(r+1)!(n+r-1)!%
}$ independent basis elements $dx^{k_{1}...k_{r+1}}$ is understood. Both ($%
i_{1},...,i_{r-1},j_{1},...,j_{n-r+1}$) and ($%
k_{1},...,k_{r+1},l_{1},...,l_{n-r-1}$) constitute even permutations of ($%
1,...,n$). In paper II, the constants $\mu _{1}$ and $\mu _{2}$ will be
determined, $r_{12}^{\lambda }$ will be specified, and proof of this theorem
will be provided.

By the linear nature of the theorem, it is clear that the theorem extends to
differential forms of inhomogeneous grade. The summations would then apply,
in addition, to all possible values of $r$. The preceding equations would
then take the still simpler form%
\begin{equation}
\alpha (x)=\mu d\left[ dx^{A}I^{\delta (A)}\right] +\mu \delta \left[
dx^{A}I^{d(A)}\right] ,
\end{equation}%
\begin{equation}
I^{\delta (A)}\equiv \int_{E_{n}^{\prime }}\frac{1}{r_{12}^{\lambda }}%
(\delta ^{\prime }\alpha ^{\prime })\wedge dx^{\prime \bar{A}},
\end{equation}%
\begin{equation}
I^{d(A)}\equiv \int_{E_{n}^{\prime }}\frac{1}{r_{12}^{\lambda }}(d^{\prime
}\alpha ^{\prime })\wedge dx^{\prime \bar{A}},
\end{equation}%
where $A$ labels the basis in the K\"{a}hler algebra, and where $dx^{\prime 
\bar{A}}$ is meant to be the element of that basis such that $dx^{\prime
A}\wedge dx^{\prime \bar{A}}$ is the basis unit element of grade $n.$ Of
course, these equations reduce to the previous ones if $\alpha $ is of
homogeneous grade.

One can then use that theorem to produce the announced major improvement on
Hodge's decomposition theorem with the help of embedding and through the use
of the KC. Retrospectively, the result obtained amounts to the integration
of the system constituted by the specification of the exterior differential
and co-differential of an arbitrary differential form on Riemannian
manifolds to which Stokes theorem can be applied (of course, one would have
to actually perform the integrations in each case, as is the case also with
Helmholtz theorem). A preprint ``Helmholtz Theorem for Differential Forms in
Arbitrary Euclidean Space, and a Hodge Theorem that Explicitly Exhibits the
Decomposition Terms'' will soon be ready for posting in arXiv math-analysis
if a reader qualified in that area volunteers to sponsor it.

What I have just stated may perhaps seem a little bit unbelievable. It
should not be so. Unbelievable is that KC has so long being overlooked. It
is a superb calculus. Using it, I have published several applications in
different fields of mathematics and physics. The last one was a paper in
high energy physics \cite{V51} which would be followed by posting of ``U(1) $%
\times $ SU(2) $\times $ SU(3) from the tangent bundle'' if I found a
sponsor for the HEP-theory section, sponsorship which I hereby request to
those qualified and who may have consulted that reference. .

\end{document}